\documentclass[10pt, a4paper]{article}
\usepackage{amscd}
\usepackage{amsmath}
\usepackage{amsfonts}
\usepackage{caption}
\begin{document}

\date{}

\newtheorem{theorem}{Theorem}
\newtheorem{corollary}[theorem]{Corollary}
\newtheorem{definition}[theorem]{Definition}
\newtheorem{lemma}[theorem]{Lemma}
\newtheorem{proposition}[theorem]{Proposition}
\newtheorem{remark}[theorem]{Remark}
\newtheorem{example}[theorem]{Example}
\newtheorem{notation}[theorem]{Notation}
\def\Qed{\hfill\raisebox{.6ex}{\framebox[2.5mm]{}}\\[.15in]}

\title{On the computation of singular plane curves\\ and quartic surfaces}

\author{Carlos Rito}

\pagestyle{myheadings}
\maketitle
\setcounter{page}{1}

\begin{abstract}

Two Magma functions are given: one computes linear systems of plane curves with non-ordinary singularities and the other computes a scheme which parametrizes given degree plane curves with given singularities. 
These functions provide an efficient tool to construct explicit equations of singular plane algebraic curves.

By computing singular branch curves, we obtain equations of normal quartic surfaces in $\mathbb C\mathbb P^3$ having the following combinations of rational double points:
$\mathsf D_5\mathsf E_7\mathsf E_7,$ $\mathsf D_7\mathsf D_6\mathsf D_6,$ $\mathsf E_6\mathsf D_8\mathsf D_5,$
$\mathsf E_6\mathsf D_{13},$ $\mathsf E_6\mathsf E_6\mathsf D_7,$ $\mathsf E_6\mathsf E_8\mathsf D_5,$
$\mathsf E_7\mathsf D_6\mathsf D_6,$ $\mathsf E_7\mathsf D_{12},$ $\mathsf E_7\mathsf E_6\mathsf D_6.$
These are all possible cases with total Milnor number $19$ which have no point of type $\mathsf A_n.$

\noindent 2000 MSC: 14Q05, 14H50.
\end{abstract}

\section{Introduction}
A quartic surface in $\mathbb C\mathbb P^3$ with only rational double points as singularities is an example of a $K3$ surface (Calabi-Yau variety of complex dimension 2), an important class of surfaces with applications to theoretical physics, in particular to string theory.
Mathematicians have been studying $K3$ surfaces for over one hundred years, at least since \cite{Hud}.

Several authors have studied normal quartic surfaces with rational double points.
A Dynkin graph of such a double point is of type $\mathsf A_n,$ $\mathsf D_m$ or $\mathsf E_p,$ $p=6, 7$ or $8,$ with Milnor number the index $n,$ $m$ or $p.$
Let $r$ be the total Milnor number, i.e. the sum of all Milnor numbers. The maximum value for $r$ is $19.$ 
Based on Nikulin's (\cite{Ni}) and Urabe's (\cite{Ur1}, \cite{Ur2}) work, Yang (\cite{Ya}) has computed, for $r=19, 18$ or $17$, all possible configurations of rational double points that occur on a normal projective quartic surface in $\mathbb C\mathbb P^3.$

The number of possibilities for $r=19$ is $278,$ but, to my knowledge, equations of such singular surfaces have not been given (Urabe's method is not constructive), except for the case of a quartic surface with an $\mathsf A_{19}$ singularity constructed by Kato and Naruki (\cite{KN}). Surprisingly, removing from Yang's list the cases which contain at least one singularity of type $\mathsf A_n,$ only the following $9$ cases remain:
$\mathsf D_5\mathsf E_7\mathsf E_7,$ $\mathsf D_7\mathsf D_6\mathsf D_6,$ $\mathsf E_6\mathsf D_8\mathsf D_5,$
$\mathsf E_6\mathsf D_{13},$ $\mathsf E_6\mathsf E_6\mathsf D_7,$ $\mathsf E_6\mathsf E_8\mathsf D_5,$
$\mathsf E_7\mathsf D_6\mathsf D_6,$ $\mathsf E_7\mathsf D_{12},$ $\mathsf E_7\mathsf E_6\mathsf D_6.$
In this paper we compute an equation for each of these cases. We notice that our method can be used to construct other surfaces in the list. For instance an equation for the case $\mathsf A_{17}\mathsf A_1\mathsf A_1$ is not difficult to obtain. 

Let $X$ be a quartic surface with a double point $p.$ The projection from $p$ gives $X$ as a double covering of the plane.
Such a covering is determined by its {\em branch locus} $B$ (the projection of the ramification curve to $\mathbb P^2$).
After computing the defining equation of $B,$ it is not difficult to obtain the equation of $X.$

In this paper we give an efficient tool to compute singular plane curves.
All computations are implemented with the Computational Algebra System MAGMA (\cite{BCP}).
The Magma function {\em LinearSystem} computes only linear systems of curves with ordinary singularities. To overcome this restriction, we define the function $LinSys$ which calculates systems of curves having (any type of) non-ordinary singularities.

Let $L$ be a linear system of plane curves of degree $d$. To obtain a curve of $L$ with given singularities, one imposes conditions to its elements. If the number of conditions is greater than the dimension of $L,$ this curve may not exist.
Suppose that the function $LinSys$ returns no sections (or returns a non-reduced curve). In this case the singular points $p_1,\ldots,p_n$ (possibly infinitely near) may be in a special position. We give a function $ParSch$ whose output is a scheme which parametrizes plane curves (or linear systems of curves) with given singularities at $p_1,\ldots,p_n$.

These two functions can be used to construct any type of singular plane curves, depending the success mostly on computer power.

The paper is organized as follows.
In section \ref{cpc} we introduce the functions $LinSys$ and $ParSch.$ The corresponding code lines, implemented with the Computational Algebra System Magma (\cite{BCP}), are given in the Appendix. In Section \ref{qs} we relate the quartic surface $X$ with its branch locus $B$ and we present a table with possibilities for the singularities of $B.$ Finally Section \ref{eq} contains the computation of the equations of the singular quartic surfaces. We give details for two cases, being the others analogous.

\bigskip
\noindent{\bf Notation}

By {\em curve} we mean algebraic curve over the complex numbers.
An $(m_1,m_2,\ldots)$-point, or point of type $(m_1,m_2,\ldots),$
is a singular point of multiplicity $m_1,$ which resolves to a point
of multiplicity $m_2$ after one blow-up, etc. A {\em tangent direction} of a plane curve singularity of multiplicity $n$ is the direction of a line which cuts the singularity with intersection number greater than $n.$

The rest of the notation is standard in Algebraic Geometry.\\

\bigskip
\noindent{\bf Acknowledgements}

The author wishes to thank Margarida Mendes Lopes for all the support.
He is a member of the Mathematics Center of the Universidade de Tr\'as-os-Montes e Alto Douro and is a collaborator of the Center for Mathematical Analysis, Geometry and Dynamical Systems of Instituto Superior T\'ecnico, Universidade T\' ecnica de Lisboa.
This research was partially supported by FCT (Portugal) through Project PTDC/MAT/099275/2008.

\section{Computation of plane curves}\label{cpc}

\subsection{Function $LinSys$}\label{LinSys}

Let $L$ be a linear system of plane curves.
The Magma function {\em LinearSystem} computes subsystems of $L,$ but only of curves with ordinary singularities. To overcome this restriction we define the function $LinSys(L,p,m,t)$ (see the Appendix), which calculates the subsystem containing the sections which have a point of type $m[i]=(m[i]_1,\ldots,m[i]_j)$ at $p[i]$ with tangent directions given by the sequence $t[i]$ of vectors, $i=1,\ldots,\#p.$
Basically this function computes the necessary blow-ups and uses the Magma function $LinearSystem$ after each one. Then it blow-downs to return to $\mathbb P^2.$

In this function the blow-up at a point with coordinates $(x,y)=(a,b)$ is given by evaluating the elements of $L$ at $(x,(x-a)y+b),$ except if the tangent direction of the singularity is the one of the vector $(0,1)$ when it is given by evaluating at $((y-b)x+a,y).$
So, by default, the exceptional divisor corresponding to the blow-up at $(a,b)$ is the line of equation $x=a,$ except if the tangent direction is $(0,1),$ when it is the line of equation $y=b.$

For example, a linear system of cubic curves with an assigned cusp can be obtained as follows:
\begin{verbatim}
A<x,y>:=AffineSpace(Rationals(),2);
LinSys(LinearSystem(A,3),A![0,0],[2,1,1],[[1,1],[0,1]]);
\end{verbatim}

\begin{figure}[h]

$$
\begin{array}{lll}

\begin{picture}(45,25)
\qbezier(-3,37)(0,40)(0,40)
\qbezier(3,37)(0,40)(0,40)
\qbezier(37,-3)(40,0)(40,0)
\qbezier(37,3)(40,0)(40,0)
\qbezier(0,0)(10,10)(8,28)
\qbezier(0,0)(10,10)(28,8)
\qbezier(0,0)(20,0)(40,0)
\qbezier(0,0)(0,20)(0,40)
\end{picture}

\hspace{1cm}

&

\begin{picture}(45,25)
\qbezier(-3,37)(0,40)(0,40)
\qbezier(3,37)(0,40)(0,40)
\qbezier(37,-3)(40,0)(40,0)
\qbezier(37,3)(40,0)(40,0)
\qbezier(0,0)(20,0)(40,0)
\multiput(0,0)(0,4.25){10}{\line(0,1){2}}
\qbezier(10,7)(-10,20)(10,33)
\put(-7,33){\makebox(0,0){\tiny $-1$}}
\end{picture}

\hspace{1cm}

&

\begin{picture}(45,30)
\qbezier(-3,37)(0,40)(0,40)
\qbezier(3,37)(0,40)(0,40)
\qbezier(37,-3)(40,0)(40,0)
\qbezier(37,3)(40,0)(40,0)
\multiput(0,0)(0,4.25){10}{\line(0,1){2}}
\multiput(0,20)(4.25,0){10}{\line(1,0){2}}
\qbezier(0,0)(20,0)(40,0)
\qbezier(-12,8)(0,20)(12,32)
\put(-7,33){\makebox(0,0){\tiny $-2$}}
\put(35,15){\makebox(0,0){\tiny $-1$}}
\end{picture}

\end{array}
$$

\caption*{Resolution of the cusp}
\end{figure}
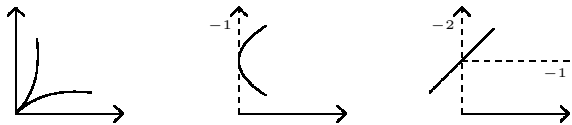

\subsection{Procedure $CndMt$}\label{CndMt}

Let $J$ be a linear system of plane curves. We want to compute the set of tuples $(p_1,\ldots,p_n)$ such that there is at least an element of $J$ with given singularities at the points $p_1,\ldots,p_n.$ Thus we have to impose conditions, on generic points, to the elements of $J.$ These conditions are minors of a matrix.

For example, suppose we want to find the elements of $J$ which have a tacnode (i.e. a $(2,2)$-point).
Let $(a,b), (u,v)$ be the coordinates of points $p_1, p_2$ such that $p_2$ is infinitely near to $p_1$. We suppose that the chart of the blow-up at $p_1$ is such that the corresponding exceptional curve has equation $x-a$. Notice that then $a=u.$ Denote the equations of the sections of $J$ by $j_1,\ldots,j_m$ and let $j_i^0$ be the equation obtained from $j_i$ by blowing-up at $p_1,$ i.e. by evaluating $j_i$ at $(x,(x-a)y+b).$ There are unique polynomials $h_i, g_i, k_i$ such that $j_i^0=h_i(x-a)^2+g_i(x-a)+k_i$ and $g_i, k_i$ are not divisible by $x-a.$ If the curves of $J$ have a double point at $p_1,$ then $k_i=g_i=0,$  $\forall i.$ To obtain double points at $p_1, p_2$ we impose conditions, vanishing of derivatives, to the equations $j_i$ and $h_i,$ $i=1,\ldots,m.$ Consider the matrix
$$
{\rm Mt}:=
\left[
\begin{array}{ccc}
j_1(a,b) & \ldots & j_m(a,b)\\
\frac{\partial j_1}{\partial x}(a,b) & \ldots & \frac{\partial j_m}{\partial x}(a,b)\\
\frac{\partial j_1}{\partial y}(a,b) & \ldots & \frac{\partial j_m}{\partial y}(a,b)\\
h_1(u,v) & \ldots & h_m(u,v)\\
\frac{\partial h_1}{\partial x}(u,v) & \ldots & \frac{\partial h_m}{\partial x}(u,v)\\
\frac{\partial h_1}{\partial y}(u,v) & \ldots & \frac{\partial h_m}{\partial y}(u,v)
\end{array}
\right]
.$$
If there is a tacnode at $p_1,$ the maximal minors of this matrix vanish. 

To define such a matrix of conditions in a more general case, involving curves with more complicated singularities, we give in the Appendix the Magma procedure
$CndMt(R,J,P,M,T,\verb+~+E,\verb+~+Mt).$

Here $R$ is a polynomial ring, $J$ is a linear system of plane curves, $P$ is a sequence of points, an element of $M$ is a sequence of multiplicities and an element of $T$ is a sequence of vectors, corresponding to tangent directions of infinitely near singularities at an element of $P.$

As output we have the matrix $Mt$ and the sequence of polynomials $E.$
In the function given in Section \ref{ParSch} below, we obtain infinitely near points by imposing $E=0.$

\subsection{Function $ParSch$}\label{ParSch}

The variety $S$ of tuples $(p_1,\ldots,p_n)$ such that there is an element of $J$ (possibly non-reduced) with given singularities at $p_1,\ldots,p_n$ is defined by:
\begin{description}
\item[$\cdot$] the vanishing of the maximal minors of $Mt$;
\item[$\cdot$] $E=0$ (infinitely near points).
\end{description}

If we are able to compute a point in $S,$ then we use the function $LinSys,$ described in Section \ref{LinSys}, to compute the sections of $J$ which have the given singularities at $p_1,\ldots,p_n.$

In the Appendix we define the Magma function $ParSch(J,P,M,T,Eq,Ne,d)$ which computes such a scheme $S.$ Here $J$ is a sequence of linear systems and the number $d$ means that the last $d$ linear systems of $J$ do not have certain singularities at the points in $P.$ In this case some of the defining equations of $S$ are of the type
$$
(1+n\alpha_1)\cdots(1+n\alpha_d)=0,
$$
where each $\alpha_i$ is a maximal minor of a certain matrix.
This is useful to obtain reduced curves.

The input $Ne$ is used to obtain no repetition on the points. For example, suppose we have points $p_1,\ldots,p_4$ (possibly infinitely near) and we want $p_1\ne p_3$ and $p_2\ne p_4.$ Then we use $Ne=[[1,3],[2,4]].$ Analogously $Eq$ is used to obtain repetition.

As in the function $LinSys,$ the exceptional divisor corresponding to the blow-up at $(a,b)$ is the line of equation $x=a,$ except if the tangent direction is $(0,1),$ when it is the line of equation $y=b.$

\section{Quartic surfaces}\label{qs}

\subsection{Quartic double planes}

Let $X\subset\mathbb P^3$ be a normal quartic surface with rational double points as its only singularities. Choosing coordinates $(x,y,z,w)$ such that $X$ has a double point $p$ at $[0:0:0:1]$, the surface is given by an homogeneous equation $$a_2w^2+b_3w+c_4=0,$$ where $a_2, b_3, c_4$ are polynomials in $(x,y,z).$ The projection of $X$ from $p$ to $\mathbb P^2$ is a $2:1$ map, ramified over the sextic curve $B$ with equation $b_3^2-4a_2c_4=0.$
The image of $p$ in $\mathbb P^2$ is the conic $C$ with equation $a_2=0,$ which is tangent to the curve $B.$ 

Let $X'$ be the double cover of $\mathbb P^2$ with branch locus $B,$ i.e. the surface given by $t^2=b_3^2-4a_2c_4$ in the weighted projective space $\mathbb P[3,1,1,1],$ with coordinates $(t,x,y,z).$
The pullback of each irreducible component of $C$ in $X'$ is the union of two $(-2)$-curves. One can obtain $X$ contracting one of these $(-2)$-curves, for each component.

Thus, given such curves $B$ and $C,$ to obtain the equation of $X$ it suffices to factor the defining polynomial of $B$ has $b_3^2-4a_2c_4.$ This is not a difficult task. A Magma function which computes this factorization and gives the equation of $X$ can be found at
\verb+http://home.utad.pt/~crito/+ .

\subsection{Singularities of double coverings}

Each singularity $q$ of $X$ is an $\mathsf A$-$\mathsf D$-$\mathsf E$ surface singularity.
The corresponding singularity of the branch locus $B$ is an $\mathsf A$-$\mathsf D$-$\mathsf E$ singularity of the same type,
except if $q$ is the projection point $p$ (\cite[III. 7.]{BHPV}).
In this case the branch curve $B$ is tangent to the conic with equation $a_2=0.$
A generic projection gives rise to a singularity in $B$ as in Table \ref{xxx} (column $[q=p]$) and
conversely one can verify that the canonical resolution of $B$ gives graphs as indicated in the table.
Here the data $m$ (multiplicities) and $t$ (tangent directions) are input data for the Magma functions $LinSys$ and $ParSch.$

\begin{table}

$$
\begin{array}{|l|l|l|}

\hline

&

q\ne p

&

q=p

\\ & &\\
\hline

\begin{array}{c}
\mathsf A_1
\end{array}

&

\begin{picture}(45,25)
\qbezier(10,-10)(20,0)(30,10)
\qbezier(10,10)(20,0)(30,-10)
\end{picture}

&

\begin{picture}(45,25)
\qbezier(2,9)(7,1)(11,9)
\qbezier(15,9)(20,1)(24,9)
\qbezier(28,9)(33,1)(37,9)
\qbezier(2,-9)(7,-1)(11,-9)
\qbezier(15,-9)(20,-1)(24,-9)
\qbezier(28,-9)(33,-1)(37,-9)
\qbezier[25](3,5)(20,5)(37,5)
\qbezier[25](3,-5)(20,-5)(37,-5)
\qbezier[10](3,-5)(-5,0)(3,5)
\qbezier[10](37,-5)(45,0)(37,5)
\end{picture}

\\ & &\\
\hline

\begin{array}{c}
\mathsf A_{2n-1}\\
n\geq 2
\end{array}

&

\begin{picture}(45,25)
\qbezier(5,10)(20,-10)(35,10)
\qbezier(5,-10)(20,10)(35,-10)
\end{picture}

m=\underbrace{[2,\ldots,2]}_{n\ {\rm times}}

&

\begin{picture}(45,30)
\qbezier(5,10)(20,-10)(35,10)
\qbezier(5,-10)(20,10)(35,-10)
\qbezier(25,10)(35,10)(35,20)
\qbezier(25,-10)(35,-10)(35,-20)
\qbezier(15,10)(5,10)(5,20)
\qbezier(15,-10)(5,-10)(5,-20)
\qbezier[30](0,-20)(20,0)(40,20)
\qbezier[30](0,20)(20,0)(40,-20)

\end{picture}

m=\underbrace{[2,\ldots,2]}_{n-1\ {\rm times}}

\\ & &\\
\hline

\begin{array}{c}
\mathsf A_{2}
\end{array}

&

\begin{picture}(45,25)
\qbezier(20,0)(30,0)(35,10)
\qbezier(20,0)(30,0)(35,-10)
\end{picture}

\begin{array}{c}
m=[2,1,1]\\
t=[[\ ],[0,1]]
\end{array}

&

\begin{picture}(45,25)
\qbezier(2,9)(7,1)(11,9)
\qbezier(15,9)(20,1)(24,9)
\qbezier(28,9)(33,1)(37,9)
\qbezier(2,-9)(7,-1)(11,-9)
\qbezier(15,-9)(20,-1)(24,-9)
\qbezier(28,-9)(33,-1)(37,-9)
\qbezier[25](0,5)(20,5)(40,5)
\qbezier[25](0,-5)(20,-5)(40,-5)
\end{picture}

\\ & & \\
\hline

\begin{array}{l}
\mathsf A_{2n}\\
n\geq 2
\end{array}

&

\begin{picture}(45,25)
\qbezier(20,0)(30,0)(35,10)
\qbezier(20,0)(30,0)(35,-10)
\end{picture}

\begin{array}{c}
m=[\underbrace{2,\ldots,2}_{n\ {\rm times}},1,1]\\
t=[[\ ],\ldots,[0,1]]
\end{array}

&

\begin{picture}(45,25)
\qbezier(25,10)(35,10)(35,20)
\qbezier(25,-10)(35,-10)(35,-20)
\qbezier(15,10)(5,10)(5,20)
\qbezier(15,-10)(5,-10)(5,-20)
\qbezier[30](0,-20)(20,0)(40,20)
\qbezier[30](0,20)(20,0)(40,-20)
\qbezier(20,0)(30,0)(35,10)
\qbezier(20,0)(30,0)(35,-10)
\end{picture}

\begin{array}{c}
m=[\underbrace{2,\ldots,2}_{n-1\ {\rm times}},1,1]\\
t=[[\ ],\ldots,[0,1]]
\end{array}

\\ & & \\
\hline

\begin{array}{c}
\mathsf D_{2n+4}\\
n\geq 0
\end{array}

&

\begin{picture}(45,25)
\qbezier(5,10)(20,-10)(35,10)
\qbezier(5,-10)(20,10)(35,-10)
\put(20,-15){\line(0,1){30}}
\end{picture}

\begin{array}{c}
m=[3,\underbrace{2,\ldots,2}_{n\ {\rm times}}]
\end{array}

&

\begin{picture}(45,25)
\qbezier(0,5)(5,0)(10,-5)
\qbezier(0,-5)(5,0)(10,5)
\qbezier(20,8)(30,-8)(40,8)
\qbezier(22,10)(38,0)(22,-10)
\qbezier(38,10)(22,0)(38,-10)
\qbezier[25](0,0)(20,0)(40,0)
\end{picture}

m=[3,\underbrace{2,\ldots,2}_{n-1\ {\rm times}}]

\\ & & \\
\hline

\begin{array}{c}
\mathsf D_{5}
\end{array}

&

\begin{picture}(45,25)
\qbezier(20,0)(30,0)(35,10)
\qbezier(20,0)(30,0)(35,-10)
\put(20,-15){\line(0,1){30}}
\end{picture}

\begin{array}{c}
m=[3,1,1]\\
t=[[\ ],[0,1]]
\end{array}

&

\begin{picture}(45,25)
\qbezier(0,5)(5,0)(10,-5)
\qbezier(0,-5)(5,0)(10,5)
\qbezier(20,8)(30,-8)(40,8)
\qbezier(20,-8)(30,8)(40,-8)
\qbezier[25](0,0)(20,0)(40,0)
\end{picture}

m=[2,2]

\\ & & \\
\hline

\begin{array}{c}
\mathsf D_{2n+5}\\
n\geq 1
\end{array}

&

\begin{picture}(45,25)
\qbezier(20,0)(30,0)(35,10)
\qbezier(20,0)(30,0)(35,-10)
\put(20,-15){\line(0,1){30}}
\end{picture}

\begin{array}{c}
m=[3,\underbrace{2,\ldots,2}_{n\ {\rm times}},1,1]\\
t=[[\ ],\ldots,[0,1]]
\end{array}

&

\begin{picture}(45,25)
\qbezier(0,5)(5,0)(10,-5)
\qbezier(0,-5)(5,0)(10,5)
\qbezier(20,8)(30,-8)(40,8)
\qbezier(30,0)(30,-4)(40,-8)
\qbezier(30,0)(30,-4)(20,-8)
\qbezier[25](0,0)(20,0)(40,0)
\end{picture}

\begin{array}{c}
m=[3,\underbrace{2,\ldots,2}_{n-1\ {\rm times}},1,1]\\
t=[[\ ],\ldots,[0,1]]
\end{array}

\\ & & \\
\hline

\begin{array}{c}
\mathsf E_6
\end{array}

&

\begin{picture}(45,25)
\qbezier(5,10)(20,-10)(35,10)
\multiput(0,0)(4.25,0){10}{\line(1,0){2}}
\put(7,-6){\makebox(0,0){\tiny $-1$}}
\end{picture}

\begin{array}{c}
m=[3,1,1,1]\\
t=[[\ ],[0,1],[0,1]]
\end{array}

&

\begin{picture}(45,25)
\qbezier(5,8)(20,-8)(35,8)
\qbezier(5,-8)(20,8)(35,-8)
\qbezier[25](0,0)(20,0)(40,0)
\put(20,-15){\makebox(0,0){\tiny ${\rm i.n.}=3+3$}}
\end{picture}

m=[2,2,2]

\\ & & \\
\hline

\begin{array}{c}
\mathsf E_7
\end{array}

&

\begin{picture}(45,25)
\qbezier(20,0)(30,0)(35,10)
\qbezier(20,0)(30,0)(35,-10)
\put(0,0){\line(1,0){40}}
\end{picture}

\begin{array}{c}
m=[3,2,1]\\
t=[[\ ],[0,1]]
\end{array}

&

\begin{picture}(45,25)
\qbezier(2,8)(20,-8)(38,8)
\qbezier(8,-8)(20,8)(32,-8)
\qbezier[25](0,0)(20,0)(40,0)
\put(20,-10){\line(0,1){20}}
\put(20,-15){\makebox(0,0){\tiny ${\rm i.n.}=1+2+3$}}
\end{picture}

m=[3,2,1]

\\ & & \\
\hline

\begin{array}{c}
\mathsf E_8
\end{array}

&

\begin{picture}(45,25)
\qbezier(20,0)(30,0)(40,10)
\qbezier(20,0)(30,0)(40,-10)
\multiput(0,0)(4.25,0){10}{\line(1,0){2}}
\put(7,-6){\makebox(0,0){\tiny $-1$}}
\end{picture}

\begin{array}{c}
m=[3,2,1,1]\\
t=[[\ ],[0,1],[1,0]]
\end{array}

&

\begin{picture}(45,25)
\qbezier(20,0)(30,0)(35,10)
\qbezier(20,0)(30,0)(35,-10)
\qbezier[25](0,0)(20,0)(40,0)
\qbezier(2,-8)(20,8)(38,-8)
\put(20,-15){\makebox(0,0){\tiny ${\rm i.n.}=3+3$}}
\end{picture}

\begin{array}{c}
m=[3,2,1]\\
t=[[\ ],[0,1]]\\
t=[[\ ],[1,0]]
\end{array}

\\ & & \\
\hline

\end{array}
$$

\caption{
Singularities of branch curve
}
\caption*{
\begin{picture}(25,5)
\multiput(0,2.5)(4.25,0){5}{\line(1,0){2}}
\end{picture}
exceptional line with self-intersection $-1$
}
\vspace{-.2cm}
\caption*{
\begin{picture}(25,5)
\qbezier[12](0,2.5)(10,2.5)(20,2.5)
\end{picture}
the image of $\{a_2=0\}$ in $\mathbb P^2$
}
\vspace{-.2cm}
\caption*{
i.n. $=$ 'intersection number'
}
\label{xxx}
\end{table}

\section{Equations}\label{eq}

In this section we compute equations of singular quartic surfaces for the cases referred in the Introduction. This is done by first computing the equation of the sextic branch curve $B$. In each case the first singularity corresponds to the projection point. For example, in a case $\mathsf E_n\mathsf E_i\mathsf D_j$ we compute a sextic curve with singular points as in Table \ref{xxx}, row $\mathsf E_n,$ column $[p=q]$ and rows $\mathsf E_i,$ $\mathsf D_j,$ column $[p\ne q].$

\subsection{$\mathsf D_5\mathsf E_7\mathsf E_7$}\label{sec1}

In this case $B=B_1+B_2$ is the union of two cubics containing points $p_1,\ldots,p_4$ such that:\\
$B_1,$ $B_2$ are tangent to the line $p_1p_2$ at $p_1;$\\
$B_1$ has a cusp at $p_3$ tangent to $B_2;$\\
$B_2$ has a cusp at $p_4$ tangent to $B_1;$\\
$p_1,$ $p_2$ are not in the line $p_3p_4$\\
(without this last condition we would obtain a non-reduced curve).

We want to find $9$ points $p_1, p_1', p_2, p_3, p_3', p_3''$ and $p_4, p_4', p_4''$ with $p_i'$ infinitely near to $p_i$ and $p_i''$ infinitely near to $p_i'.$ We fix all points except $p_4$ and $p_4'.$ Notice that each infinitely near point is defined by the tangent direction of the corresponding singularity.
%
\begin{verbatim}
A<x,y>:=AffineSpace(Rationals(),2);
L:=[LinearSystem(A,3),LinearSystem(A,3),LinearSystem(A,1), \
LinearSystem(A,1)];
P:=[A![0,0],A![1,0],A![2,1]];
M:=[
    [[1,1],[1],[2,1,1],[1,1,0]],\
    [[1,1],[1],[1,1,0],[2,1,1]],\
    [[0,0],[1],[1,0,0],[1,0,0]],\
    [[1,0],[0],[1,0,0],[1,0,0]]];
T:=[[[1,0]],[],[[1,2],[0,1]],[[],[0,1]]];
S:=ParSch(L,P,M,T,[],[[1,3,4,7]],2);
\end{verbatim} 
With the input [1,3,4,7] we obtain no repetition among the points in positions 1, 3, 4 and 7 (i.e. $p_1,\ldots,p_4$), so we get $p_4\ne p_1, p_2, p_3.$ The previous sequence [ ] is empty because we do not want any points to be equal.

Now it only remains to find a solution in $S.$
Since $S$ has dimension $1,$ we define a zero-dimensional subscheme $S_1$ and compute its points.
\begin{verbatim}
R:=Ambient(S);
S1:=Scheme(S,[R.2-2]);
PointsOverSplittingField(S1);
\end{verbatim} 
The first four coordinates of the ambient space $R$ of $S$ correspond to $p_4$ and $p_4'.$
We choose a solution and compute the curve $B$ with the given singularities at $p_1,\ldots,p_4.$
\begin{verbatim}
P:=[A![0,0],A![1,0],A![2,1],A![22/7,2]];
M:=[[2,2],[2],[3,2,1],[3,2,1]];
T:=[[[1,0]],[],[[1,2],[0,1]],[[1,28/23],[0,1]]];
J:=LinSys(LinearSystem(A,6),P,M,T);
\end{verbatim} 
This linear system $J$ has only one section, which factors as
\begin{verbatim}
(x^3 - 323/63*x^2*y - x^2 + 512/63*x*y^2 +
92/21*x*y - 254/63*y^3 - 88/21*y^2 - 2/7*y)
(x^3 - 515/126*x^2*y - x^2 + 2482/441*x*y^2 +
58/21*x*y - 1129/441*y^3 - 317/147*y^2 + 2/7*y)
\end{verbatim}
From this branch curve we obtain the equation of the quartic surface $X$ with a $\mathsf D_5$ point and two $\mathsf E_7$ points:
\begin{verbatim}
w^2*y^2 + w*x^3 - 129/28*w*x^2*y - w*x^2*z + 25/7*w*x*y*z
- 856343/254016*x^4 + 1907707/111132*x^3*y - 348881/18522*x^2*y^2
+ 302119/27783*x*y^3 - 143383/55566*y^4 + 50963/10584*x^3*z
- 126379/6174*x^2*y*z + 48976/3087*x*y^2*z - 89935/18522*y^3*z
- 793/588*x^2*z^2 + 5224/1029*x*y*z^2 - 4433/2058*y^2*z^2
- 17/147*x*z^3 + 299/2058*y*z^3 + 1/49*z^4
\end{verbatim}

\subsection{$\mathsf E_6\mathsf D_{13}$}

Here $B$ is an irreducible sextic curve with a $(2,2,2)$-point and a $(3,2,2,2,2)$-point which resolves to a cusp after $4$ blow-ups. To obtain a reduced curve we impose that $B$ is not a cubic with multiplicity $2$ nor a conic with multiplicity $3.$
\begin{verbatim}
A<x,y>:=AffineSpace(Rationals(),2);
L:=[LinearSystem(A,6),LinearSystem(A,3),LinearSystem(A,2)];
P:=[A![0,0],A![0,1]];
M:=[[[2,2,2],[3,2,2,2,2,1,1]],\
    [[1,1,1],[2,1,1,1,1,0,0]],\
    [[1,1,0],[1,1,1,1,1,1,0]]];
T:=[[[1,0],[1,0]],[[1,0],[1,1],[],[],[],[0,1]]];
S:=ParSch(L,P,M,T,[],[],2);
\end{verbatim} 
Now we proceed as in Section \ref{sec1}.
The equation of the quartic surface $X$ with an $\mathsf E_6$ point and a $\mathsf D_{13}$ point is
\begin{verbatim}
w^2*y^2 + w*x^3 - 9/256*w*x^2*y + 69/64*w*x*y*z + 111/64*w*y*z^2
+ 141393/262144*x^4 - 231/256*x^3*y - 33/128*x^2*y^2
+ 489/512*x*y^3 - 25/32*y^4 + 57363/32768*x^3*z + 33/64*x^2*y*z
- 729/256*x*y^2*z + 99/32*y^3*z + 75/32768*x^2*z^2
+ 1449/512*x*y*z^2 - 147/32*y^2*z^2 - 21/8192*x*z^3 + 97/32*y*z^3
+ 33/16384*z^4
\end{verbatim}

\subsection{The remaining equations}
Now we give the remaining equations of the quartic surfaces $X$ with rational double points with total Milnor number $19$ which have no point of type $\mathsf A_n.$ The details can be found at
\verb+http://home.utad.pt/~crito/+ .
\\\\
$\mathsf D_7\mathsf D_6\mathsf D_6$
\begin{verbatim}
w^2*y^2 + w*x^3 - 11/4*w*x^2*y - w*x^2*z + 4*w*x*y*z - 3*w*y*z^2
+ 113/64*x^4 - 8*x^3*y + 16*x^2*y^2 - 7/2*x*y^3 - 11/2*y^4
- 6*x^3*z + 25/2*x^2*y*z - 30*x*y^2*z + 23/2*y^3*z + 89/8*x^2*z^2
- 11/2*x*y*z^2 + 31/2*y^2*z^2 - 9*x*z^3 - 3/2*y*z^3 + 9/4*z^4
\end{verbatim}
$\mathsf E_6\mathsf D_8\mathsf D_5$
\begin{verbatim}
w^2*y^2 + w*x^3 - 3*w*x^2*y - 10/3*w*x*y*z + 7/3*w*y*z^2
- 29/12*x^4 + 15/2*x^3*y - 17/3*x^2*y^2 - 4/3*x*y^3 + 4*y^4
+ 2/3*x^3*z + 14/3*x^2*y*z + 20/3*x*y^2*z - 32/3*y^3*z
+ 5/18*x^2*z^2 - 28/3*x*y*z^2 + 8*y^2*z^2 + 1/9*x*z^3 + 1/36*z^4
\end{verbatim}
$\mathsf E_6\mathsf E_6\mathsf D_7$\ \ \ \ \ \ \verb#( r^2 - 33/73*r + 9/292 = 0 )#
\begin{verbatim}
w^2*y^2 + w*x^3 - 9/8*w*x^2*y + 3/2*w*x*y*z - 9/2*w*y*z^2
+ 1/256*(-292*r + 207)*x^4 + 1/64*(-146*r + 69)*x^3*y
- 15/32*x^2*y^2 + 1/64*(146*r - 33)*x*y^3 + 1/128*(146*r - 83)
*y^4 + 1/16*(73*r - 66)*x^3*z + 15/8*x^2*y*z + 1/32*(-438*r + 99)
*x*y^2*z + 1/16*(-146*r + 75)*y^3*z + 39/32*x^2*z^2 + 1/16
*(438*r - 99)*x*y*z^2 + 1/16*(438*r - 201)*y^2*z^2 + 1/4*(-73*r
+ 3)*x*z^3 + 1/4*(-146*r + 59)*y*z^3 + 1/16*(292*r - 21)*z^4
\end{verbatim}
$\mathsf E_6\mathsf E_8\mathsf D_5$
\begin{verbatim}
w^2*y^2 + w*x^3 - 27/16*w*x^2*y - 21/8*w*x*y*z - 27/16*w*y*z^2
- 615/1024*x^4 + 27/32*x^3*y - 39/64*x^2*y^2 + 27/32*x*y^3
- 25/64*y^4 + 135/256*x^3*z + 39/32*x^2*y*z - 81/32*x*y^2*z
+ 51/32*y^3*z + 1299/512*x^2*z^2 + 81/32*x*y*z^2
- 39/16*y^2*z^2 + 351/256*x*z^3 + 53/32*y*z^3 + 297/1024*z^4
\end{verbatim}
$\mathsf E_7\mathsf D_6\mathsf D_6$
\begin{verbatim}
w^2*y^2 + w*x^3 - 1/4*w*x^2*y + 3/4*w*x*y*z + 1/4*w*y*z^2
+ 25/64*x^4 - 5/16*x^3*y - 1/16*x^2*y^2 + 1/8*x*y^3 - 1/32*y^4
+ 11/32*x^3*z + 1/8*x^2*y*z - 5/16*x*y^2*z + 3/32*y^3*z
+ 3/64*x^2*z^2 + 1/4*x*y*z^2 - 3/32*y^2*z^2 + 1/32*x*z^3
+ 1/32*y*z^3 + 1/64*z^4
\end{verbatim}
$\mathsf E_7\mathsf D_{12}$
\begin{verbatim}
w^2*y^2 + w*x^3 - 1/2*w*x^2*y + 3/2*w*x*y*z + 1/2*w*y*z^2
+ 13/16*x^4 - 3/4*x^3*y - 1/4*x^2*y^2 + 3/4*x*y^3 - 1/4*y^4
+ 5/8*x^3*z + 1/2*x^2*y*z - 2*x*y^2*z + 3/4*y^3*z + 3/16*x^2*z^2
+ 7/4*x*y*z^2 - 3/4*y^2*z^2 - 1/8*x*z^3 + 1/4*y*z^3 + 1/16*z^4
\end{verbatim}
$\mathsf E_7\mathsf E_6\mathsf D_6$
\begin{verbatim}
w^2*y^2 + w*x^3 - 3/4*w*x^2*y - 3/4*w*x*y*z - 3/4*w*y*z^2
- 23/64*x^4 + 9/16*x^3*y - 3/16*x^2*y^2 - 1/16*x*y^3 + 1/32*y^4
- 21/32*x^3*z + 3/8*x^2*y*z + 3/16*x*y^2*z - 3/32*y^3*z
+ 15/64*x^2*z^2 - 3/16*x*y*z^2 + 3/32*y^2*z^2 + 11/32*x*z^3
- 1/32*y*z^3 + 9/64*z^4
\end{verbatim}


\appendix
\section*{Appendix: Magma code}

\footnotesize
\begin{verbatim}
function LinSys(L,p,m,t)

  A:=Ambient(L);
  J:=LinearSystem(L,p,[m[i][1]:i in [1..#m]]);

  for g in [1..#p] do
    p0:=[BaseField(A)|p[g][1],p[g][2]];
    b:=[];i:=[];

    for j:=1 to #t[g] do
      if #Sections(J) eq 0 then break;end if;
      if t[g][j][1] eq 0 then k:=[2,1];else k:=[1,2];end if;
      b:=b cat [p0];i:=i cat [k];
      Bup:=[Evaluate(Sections(J)[i],A.k[2],(A.k[1]-p0[k[1]])*A.k[2]+ \
      p0[k[2]]) div (A.k[1]-p0[k[1]])^m[g][j]:i in [1..#Sections(J)]];
      p0[k[2]]:=t[g][j][k[2]]/t[g][j][k[1]];
      J:=LinearSystem(LinearSystem(A,Bup),A!p0,m[g][j+1]);
    end for;

    for j:=#b to 1 by -1 do
      if #Sections(J) eq 0 then break;end if;
      i1:=i[j][1];i2:=i[j][2];
      Bdn:=[Evaluate((A.i1-b[j][i1])^m[g][j]*Sections(J)[i],A.i2, \
      (A.i2-b[j][i2])/(A.i1-b[j][i1])):i in [1..#Sections(J)]];   
      R:=Universe(Bdn);
      h:=hom<R->CoordinateRing(A)|[A.1,A.2]>;
      J:=LinearSystem(A,[h(Bdn[i]):i in [1..#Bdn]]);
    end for;

  end for;

  return J;

end function;
\end{verbatim} 

\begin{verbatim}
procedure CndMt(R,J,P,M,T,~E,~Mt)

  function D(F,i,j,a,b);
    P:=Parent(F);
    for n in [1..a] do F:=Derivative(F,P.i);end for;
    for n in [1..b] do F:=Derivative(F,P.j);end for;
    return F;
  end function;

  k:=0;a:=0;E:=[];Mt:=[[]:i in [1..#Sections(J)]];
  A:=Ambient(J);
  PP:=[[P[i]]:i in [1..#P]] cat [[]:i in [1..#M-#P]];

  for w in [1..#PP] do
    p:=PP[w];m:=M[w];t:=T[w];t1:=[];

    for i:=1 to #t do
      if t[i] eq [] then t1:=[t[j]:j in [i..#t]];t:=[t[j]:j in [1..i-1]];
      break;end if;
    end for;

    if #m le #t+#p then continue;end if;
    Bup:=Sections(J);
    if #p gt 0 and #m gt #t+1 then
      pa:=[BaseField(A)|p[1][1],p[1][2]];
      for j:=1 to #t do
        if t[j][1] eq 0 then q:=[2,1];else q:=[1,2];end if;
        Bup:=[Evaluate(Bup[i],A.q[2],(A.q[1]-pa[q[1]])*A.q[2]+pa[q[2]]) \
        div (A.q[1]-pa[q[1]])^m[j]:i in [1..#Bup]];
        pa[q[2]]:=t[j][q[2]]/t[j][q[1]];
      end for;
      Bup:=[Evaluate(Bup[i],A.2,(A.1-pa[1])*A.2+pa[2]) \
      div (A.1-pa[1])^m[#t+1]:i in [1..#Bup]];
    end if;

    a:=a+k;k:=2*(#m-#t-#p);
    h:=hom<PolynomialRing(J)->R|[R.(a+k-1),R.(a+k)]>;
    H:=[h(Bup)];
    if p eq [] then u:=R.(a+1);else u:=p[1][1];end if;
    su:=[u,R.(a+2)];

    o:=0;
    for j in [1..#m-#t-#p-1] do
      l:=H[j];
      H[j]:=[Evaluate(Evaluate(l[i],R.(a+k-1),R.(a+2*j-1)),R.(a+k), \
      R.(a+2*j)):i in [1..#l]];
      if t1 ne [] and t1[j+#p] ne []
        then
        if t1[j+#p][o+1] eq 0
          then o:=Abs(o-1);u:=Abs(o-1)*R.(a+2*j-1)+o*R.(a+2*j);
          su:=su cat [Abs(o-1)*u,o*u];
          else su:=su cat [Abs(o-1)*u+o*t1[j+#p][1]/t1[j+#p][2], \
          o*u+Abs(o-1)*t1[j+#p][2]/t1[j+#p][1]];
        end if;
        else su:=su cat [Abs(o-1)*u+o*R.(a+2*j+1),o*u+Abs(o-1)*R.(a+2*j+2)];
      end if;
      l:=[Evaluate(l[i],R.(a+k-o), \
      (R.(a+k-1+o)-R.(a+2*j-1+o))*R.(a+k-o)+R.(a+2*j-o)):i in [1..#l]];

      for c in [1..m[#t+#p+j]] do
        l:=[(l[i]-Evaluate(l[i],R.(a+k-1+o),R.(a+2*j-1+o))) \
        div (R.(a+k-1+o)-R.(a+2*j-1+o)):i in [1..#l]];
      end for;

      H[j+1]:=l;
    end for;

    mt:=[&cat[[D(H[j][i],a+2*j-1,a+2*j,b,c-b):b in [0..c], \
    c in [0..m[#p+#t+j]-1]]:j in [1..#H]]:i in [1..#H[1]]];
    mte:=[[Evaluate(mt[i][o],[R.i:i in [1..a]] cat su cat \
    [R.i:i in [a+k+1..Rank(R)]]):o in [1..#mt[1]]]:i in [1..#mt]];
    Mt:=[Mt[i] cat mte[i]:i in [1..#mte]];
    E:=E cat [R.(a+i)-su[i]:i in [1..k]];
  end for;

  Mt:=Matrix(Mt);

end procedure;
\end{verbatim} 

\begin{verbatim}
function EqNe(R,P,M,T,Eq,Ne,d)

P:=[[P[i]]:i in [1..#P]] cat [[]:i in [1..#M[1]-#P]];
U:=[];k:=1;
for j:=1 to #M[1] do
  if P[j] eq [] then p:=[R.k,R.(k+1)];k:=k+2;
  else p:=[P[j][1][1],P[j][1][2]];end if;
  U:=U cat p;
  o:=0;
  for h:=1 to #T[j] do
    if T[j][h] eq [] then U:=U cat [R.k,R.(k+1)];k:=k+2;
    else
      if T[j][h][o+1] eq 0 then o:=Abs(o-1);end if;
      U:=U cat [Abs(o-1)*p[1]+o*T[j][h][2-o]/T[j][h][1+o], \
      o*p[2]+Abs(o-1)*T[j][h][2-o]/T[j][h][1+o]];
    end if;
  end for;
end for;

eql:=[];
for i:=1 to #Eq do
  eql:=eql cat [U[2*Eq[i][m]-1]-U[2*Eq[i][m+1]-1]:m in [1..#Eq[i]-1]] \
  cat [U[2*Eq[i][m]]-U[2*Eq[i][m+1]]:m in [1..#Eq[i]-1]];
end for;

neq:=[];z:=Rank(R)-d;
for i:=1 to #Ne do
  for n:=1 to #Ne[i] do
    for m:=n+1 to #Ne[i] do
      neq:=neq cat [(1+R.z*(U[2*Ne[i][n]-1]-U[2*Ne[i][m]-1]))* \
      (1+R.z*(U[2*Ne[i][n]]-U[2*Ne[i][m]]))];
      z:=z-1;
    end for;
  end for;
end for;

S:=Scheme(AffineSpace(R),eql cat neq);

return S;

end function;
\end{verbatim} 

\begin{verbatim}
function ParSch(LL,P,MM,T,Eq,Ne,d)

  if Parent(MM) eq Parent([[1]]) then LL:=[LL];MM:=[MM];end if;
  T0:=T;
  for j:=1 to #T do
    for i:=1 to #T[j] do
      if T[j][i] eq [] then T0[j]:=[T[j][h]:h in [1..i-1]];break;end if;
    end for;
  end for;

  n:=&+[2*(#MM[1][i]-#T0[i]):i in [1..#MM[1]]]-2*#P+d;
  for i:=1 to #Ne do
    n:=n+#Ne[i]*(#Ne[i]-1)/2;
  end for;
  R:=PolynomialRing(BaseField(Ambient(LL[1])),Integers()!n);

  minors:=[[]:i in [1..#LL]];
  for z in [1..#LL] do
    L:=LL[z];M:=MM[z];
    N:=[[M[j][i]:i in [1..#T0[j]+1]]:j in [1..#P]];
    if #P ne 0 then
      J:=LinSys(L,P,N,[T0[i]:i in [1..#P]]);else J:=L;
    end if;
    CndMt(R,J,P,M,T,~E,~Mt);
    minors[z]:=Minors(Mt,#Sections(J));
  end for;

  df:=[];
  if d gt 0 then
    z:=Rank(R);
    for i:=#LL-d+1 to #LL do
      df:=df cat [&*[1+R.z*minors[i][j]:j in [1..#minors[i]]]];
      z:=z-1;
    end for;
  end if;

  S1:=EqNe(R,P,MM,T,Eq,Ne,d);
  S2:=Scheme(Ambient(S1),&cat[minors[i]:i in [1..#LL-d]] cat E cat df);
  S:=Intersection(S1,S2);

  return S;

end function;
\end{verbatim}

\normalsize

\bibliography{ref}

\bigskip
\bigskip

\noindent Carlos Rito
\\ Departamento de Matem\' atica
\\ Universidade de Tr\' as-os-Montes e Alto Douro
\\ 5001-801 Vila Real
\\ Portugal
\\\\
\noindent {\it e-mail:} crito@utad.pt

\end{document}